\documentclass{amsart}
\usepackage{mathrsfs,times,amsmath,amssymb,nicefrac,float}
\newtheorem{theorem}{Theorem}
\newtheorem{corollary}[theorem]{Corollary}
\newtheorem{lemma}[theorem]{Lemma}
\newtheorem{proposition}[theorem]{Proposition}

\newcommand\nc\newcommand

\nc\bfa{{\boldsymbol a}}\nc\bfA{{\bf A}}\nc\cA{{\mathcal A}}\nc\sA{{\mathscr A}}
\nc\bfb{{\boldsymbol b}}\nc\bfB{{\bf B}}\nc\cB{{\mathcal B}}\nc\sB{{\mathscr B}}
\nc\bfc{{\boldsymbol c}}\nc\bfC{{\bf C}}\nc\cC{{\mathcal C}}\nc\sC{{\mathscr C}}
\nc\bfd{{\boldsymbol d}}\nc\bfD{{\bf D}}\nc\cD{{\mathcal D}}\nc\sD{{\mathscr D}}
\nc\bfe{{\boldsymbol e}}\nc\bfE{{\bf E}}\nc\cE{{\mathcal E}}\nc\sE{{\mathscr E}}
\nc\bff{{\boldsymbol f}}\nc\bfF{{\bf F}}\nc\cF{{\mathcal F}}\nc\sF{{\mathscr F}}
\nc\bfg{{\boldsymbol g}}\nc\bfG{{\bf G}}\nc\cG{{\mathcal G}}\nc\sG{{\mathscr G}}
\nc\bfh{{\boldsymbol h}}\nc\bfH{{\bf H}}\nc\cH{{\mathcal H}}\nc\sH{{\mathscr H}}
\nc\bfi{{\boldsymbol i}}\nc\bfI{{\bf I}}\nc\cI{{\mathcal I}}\nc\sI{{\mathscr I}}
\nc\bfj{{\boldsymbol j}}\nc\bfJ{{\bf J}}\nc\cJ{{\mathcal J}}\nc\sJ{{\mathscr J}}
\nc\bfk{{\boldsymbol k}}\nc\bfK{{\bf K}}\nc\cK{{\mathcal K}}\nc\sK{{\mathscr K}}
\nc\bfl{{\boldsymbol l}}\nc\bfL{{\bf L}}\nc\cL{{\mathcal L}}\nc\sL{{\mathscr L}}
\nc\bfm{{\boldsymbol m}}\nc\bfM{{\bf M}}\nc\cM{{\mathcal M}}\nc\sM{{\mathscr M}}
\nc\bfn{{\boldsymbol n}}\nc\bfN{{\bf N}}\nc\cN{{\mathcal N}}\nc\sN{{\mathscr N}}
\nc\bfo{{\boldsymbol o}}\nc\bfO{{\bf O}}\nc\cO{{\mathcal O}}\nc\sO{{\mathscr O}}
\nc\bfp{{\boldsymbol p}}\nc\bfP{{\bf P}}\nc\cP{{\mathcal P}}\nc\sP{{\mathscr P}}
\nc\bfq{{\boldsymbol q}}\nc\bfQ{{\bf Q}}\nc\cQ{{\mathcal Q}}\nc\sQ{{\mathscr Q}}
\nc\bfr{{\boldsymbol r}}\nc\bfR{{\bf R}}\nc\cR{{\mathcal R}}\nc\sR{{\mathscr R}}
\nc\bfs{{\boldsymbol s}}\nc\bfS{{\bf S}}\nc\cS{{\mathcal S}}\nc\sS{{\mathscr S}}
\nc\bft{{\boldsymbol t}}\nc\bfT{{\bf T}}\nc\cT{{\mathcal T}}\nc\sT{{\mathscr T}}
\nc\bfu{{\boldsymbol u}}\nc\bfU{{\bf U}}\nc\cU{{\mathcal U}}\nc\sU{{\mathscr U}}
\nc\bfv{{\boldsymbol v}}\nc\bfV{{\bf V}}\nc\cV{{\mathcal V}}\nc\sV{{\mathscr V}}
\nc\bfw{{\boldsymbol w}}\nc\bfW{{\bf W}}\nc\cW{{\mathcal W}}\nc\sW{{\mathscr W}}
\nc\bfx{{\boldsymbol x}}\nc\bfX{{\bf Z}}\nc\cX{{\mathcal X}}\nc\sX{{\mathscr X}}
\nc\bfy{{\boldsymbol y}}\nc\bfY{{\bf Y}}\nc\cY{{\mathcal Y}}\nc\sY{{\mathscr Y}}
\nc\bfz{{\boldsymbol z}}\nc\bfZ{{\bf Z}}\nc\cZ{{\mathcal Z}}\nc\sZ{{\mathscr Z}}

\newcommand{\reals}{\mathbb{R}}

\def\rank{\qopname\relax{no}{rank}}

\newcommand\half{\nicefrac12}
\nc{\remove}[1]{}
\begin{document}

\title[Sets with few distances]	
{Bounds on sets with few distances} \thanks{{\em Date}\/: \today.\/ {\em AMS Subject Classification}: Primary 05E30, Secondary 94B65. }%

\author[A. Barg]{Alexander Barg$^\ast$}\thanks{$^\ast$
Dept. of Electrical and Computer Engineering and Institute for Systems 
Research, University of Maryland, College Park, MD 20742, 
and Institute for Problems of Information Transmission, Russian Academy of Sciences, Moscow, Russia. Email: abarg@umd.edu. Research 
supported in part by NSF grants DMS0807411, CCF0635271, and CCF0830699.}
\author[O~.R. Musin]{Oleg R. Musin$^\dag$}\thanks{$^\dag$ 
Department of Mathematics, University of Texas at Brownsville, 80 Fort Brown, Brownsville, TX 78520, USA. Email: omusin@gmail.com. Research 
supported in part by NSF grant DMS0807640 and NSA grant MSPF-08G-201.}

\begin{abstract}
We derive a new estimate of the size of finite sets of points 
in metric spaces with few distances. The following applications are considered:
\begin{itemize}
  \item we improve the 
Ray-Chaudhuri--Wilson bound of the size of uniform intersecting
families of subsets;
   \item we refine the bound of Delsarte-Goethals-Seidel on the maximum size of spherical sets
with few distances;
  \item we prove a new bound on codes with few distances
in the Hamming space, improving an earlier result of Delsarte.
\end{itemize}
We also find the size of maximal binary codes and maximal constant-weight
codes of small length with 2 and 3 distances.
\end{abstract}
\maketitle

\vspace*{-.3in}
\section{Introduction}

We consider finite collections of points in a metric space $X$ 
with distance function $d$. Following the terminology of coding theory
we call such collections codes.
We say that 
$\cC\subset X$ is an $s$-code if the set of distances $d(\bfx_1,\bfx_2)$
between any two distinct points of $\cC$ has size $s$. The subject of
this paper is estimates for the size (the number of points) of
$s$-codes. 

The study of $s$-codes in $\reals^n$ was initiated by 
Einhorn and Schoenberg \cite{ein66}. 
Delsarte \cite{del73,del73a} obtained several classical results for 
$s$-codes in finite spaces, while for the case of the unit sphere
$S^{n-1}\subset \reals^n$ the problem of bounding the size of $s$-codes
was first addressed by Delsarte, Goethals, and Seidel in \cite{del77b}.
Codes with few distances  
in finite spaces are closely related to the well-known combinatorial 
problem of bounding the size of families of sets with restricted 
intersections. Results of this kind are often called intersection
theorems in combinatorial literature. They have been a subject of 
extensive studies beginning with the work of
Ray-Chaudhuri and Wilson \cite{ray75}.
Their proofs are mostly based on two general methods, namely, the method of
linearly independent polynomials, see e.g., Alon et al. 
\cite{alo91}, Blokhuis \cite{blo93}, Babai et al. \cite{bab95},
and on Delsarte's linear programming method \cite{del73a,del77b}.

Recently an improvement of the Delsarte-Goethals-Seidel 
bound on spherical $s$-codes 
for the case $s=2$ was obtained in the second author's paper \cite{mus09a}. 
Following this result, Nozaki \cite{noz09a} proved a general 
bound on the size of spherical $s$-codes.
We continue this line of work, employing Delsarte's ideas to derive a 
general improvement of the
bound \cite{del77b} for every even $s$ as well as new estimates of the
size of $s$-codes over a finite alphabet. The latter result 
also enables us to tighten 
the Ray-Chaudhuri--Wilson bound on the size of uniform $s$-intersecting
families. Of course, both these bounds are known to be tight in general, 
so our improvements are only valid under some assumptions on the size of
the intersections.

\section{A bound on $s$-codes}\label{sect:general}
In this section we present a general bound on the size of $s$-codes 
(Theorem \ref{thm:new}). The bound is most conveniently described in the 
context of harmonic analysis. 
This approach to packings of metric spaces was introduced in 
\cite{del73,del77b} for finite spaces known as
association schemes and the sphere $S^{n-1}$ respectively.
It was generalized in \cite{kab78} to all distance transitive compact metric 
spaces. Under this approach the space $X$ is 
viewed as a homogeneous space of its isometry group $G.$ 
The space $X$ is called distance transitive if $G$ acts transitively 
on ordered pairs of points of $X$ at a given distance.
Denote by $d\alpha$ the normalized $G$-invariant measure on $X$. 
The space $L^2(X,d\alpha)$ of complex-valued square-integrable functions on $X$ 
decomposes into a finite or countably infinite
direct sum of pairwise orthogonal finite-dimensional linear spaces $V_i$ 
of functions called (generalized)
spherical harmonics. Let us fix a basis of spherical harmonics 
$(\phi_{i,1},\dots,\phi_{i,h_i})$ in the space $V_i,$ where $h_i=\dim V_i$. 
Since $X$ is distance transitive, the function 
  \begin{equation}\label{eq:addition}  
  P_i(\bfx,\bfy)=\sum_{j=1}^{h_i}\phi_{i,j}(\bfx)
\overline{\phi_{i,j}(\bfy)}
   \end{equation}
depends only on the distance $d(\bfx,\bfy).$ This expression
is called the addition formula in the theory of special functions, and it
is only this formula that we need in later derivations.
Below we use small $p$ to refer to functions obtained from the functions
(\ref{eq:addition}) once the pair of points $\bfx,\bfy$ is replaced
by the distance between them, and use $x$ to denote this distance.
In particular, $p_i(x)$ is a univariate real polynomial of degree $i.$
Without loss of generality we assume that $p_0\equiv 1.$

In the cases of interest to us, the functions $p_i$ form a family of 
classical orthogonal polynomials. Namely, consider the linear functional
$\sL(f)=\int f(x)d\mu(x),$ where $d\mu$ 
is the
measure induced by $d\alpha$ on the set of possible values of the distance
on $X$. Then
$\sL(p_ip_j)=0$ for $i\ne j,$ and
   $$
     r_i\sL(p_i^2)=1,    \quad \text{ where } r_i=\frac 1{h_i}.
   $$
As is well known (e.g., \cite[p.244]{and99}), the polynomials $p_i$ satisfy a three-term recurrence
of the form 
   \begin{equation}\label{eq:3term}
    xp_i=a_ip_{i+1}+b_i p_i+c_ip_{i-1},
   \end{equation}
where the numbers $a_i,b_i,c_i$ can be easily computed.
Given a polynomial $f(x)$ of degree $s$ we can compute its 
Fourier coefficients in the basis $\{p_i\}$ in a usual way, namely,
   \begin{equation}\label{eq:fc}
     f_i=r_i\sL(fp_i) \qquad(0\le i\le s).
   \end{equation}

Our primary examples will be the Hamming space $H_q^n=(Z_q)^n$
where $Z_q$ is the set of integers mod $q$,
the binary Johnson
space $J^{n,w}$ formed by the $n$-dimensional binary vectors with $w$ ones,
$w\leq n/2$, and the sphere $S^{n-1}.$ The distance in $H_q^n$ is defined as
$d_H(\bfx_1,\bfx_2)=|\{i: x_{1i}\ne x_{2i}\}|,$ the distance in $J^{w,n}$
is given by $d_J(\bfx_1,\bfx_2)=\half d_H(\bfx_1,\bfx_2),$ and the distance
on $S^{n-1}$ is measured as the inner product between the vectors.

To illustrate the above ideas, let us consider the Hamming case
$X=H_q^n.$
A typical isometry of $X$ is a permutation of coordinates
followed by a permutation of symbols in every coordinate, i.e., 
$G=S_q \wr S_n.$ 
An orthogonal basis of the space $V_i$ is formed of functions 
$\phi_{i,j}(\bfx)=e^{\frac{2\pi i}q(\alpha_1 x_{l_1}+\dots+\alpha_ix_{l_i})},$
where $1\le l_1<\dots<l_i\le n$ is an $i$-subset of $[n]$ and 
$\alpha_m\in Z_q\backslash0, m=1,\dots,i.$  
There are $h_i=\binom ni(q-1)^i$ linearly independent functions $\phi_{i,j}$
of this form.
Then $p_i$ is a Krawtchouk polynomial $K_i(x)$ of degree $i$ whose explicit
form can be found from (\ref{eq:addition}).
We have
    $$
    K_i(x)=\sum_{j=0}^i (-1)^j \binom xj\binom{n-x}{i-j}(q-1)^{i-j},
    $$
In particular, $K_i(0)=\binom ni(q-1)^i,$
  \begin{equation}\label{eq:Klow}
     \begin{array}{cc}{\displaystyle K_0(x)=1, \;K_1(x)=n(q-1)-qx},\\[.1in]
       \displaystyle{K_2(x)=\half\{q^2x^2-q(2qn-q-2n+2)x+(q-1)^2n(n-1)\}.}
\end{array}
  \end{equation}

For $X=J^{n,w}$ the polynomials $p_i$ form a certain family of discrete
Hahn polynomials \cite{del78b}. The Hahn polynomial of degree $i$
is given by
   $$
     Q_i(x)=\Big(\binom ni -\binom n{i-1}\Big) \sum_{j=0}^i(-1)^j\frac{{\binom ij}\binom{n+1-i}j}
         {\binom wj\binom {n-w}j}\binom xj.
   $$
\remove{In particular,
   \begin{equation}\label{eq:Klow}
     \begin{array}{cc}{\displaystyle K_0(x)=1, \;
    H_1(x)=(n-1)\Big(1-\frac {nx}{(n-w)w}\Big),\;
      K_2(x)=2x^2-2nx+\binom n2.}
     \\[.1in]
    {\displaystyle K_3(x)=-\frac43 x^3+2nx^2-\frac x3(3n^2-3n+2)+\binom n3.}\end{array}
  \end{equation} }
Finally, for $S^{n-1}$ the functions $p_i$ are given by the Gegenbauer
polynomials $G_i(t).$ The explicit form and properties of these polynomials
are well known. All the information about them
that we need is listed in Table \ref{table1} together with the 
corresponding properties of $K_i$ and $Q_i.$
There is no single reference with the proofs of these formulas although they
are mentioned in many places. The primary sources are Koekoek and Swarttouw
\cite{koe98} (or the recent book \cite{koe10}) or 
Andrews et al. \cite{and99}, Ch.6, but the normalizations there are
different from the ones used above. Hahn polynomials are also discussed
by Delsarte in \cite{del73} (without being identified as such) and \cite{del78b}.

\begin{table}[t]
  \begin{tabular}{lccc}
      $X$ &$H_q^n,\; q\ge 2$&$J^{n,w}$&$S^{n-1}$\\[1mm]\hline
   \vspace*{1mm}
   $d\mu$ &$q^{-n}\binom ni(q-1)^i$ &$\frac{\binom wi\binom{n-w}i}{\binom nw}$
         &$\frac{\Gamma(n/2)}{2\pi^{n/2}}(1-x^2)^{(n-3)/2}dx$\\
   $h_i$  &$\binom ni(q-1)^i$& $\binom ni-\binom n{i-1}$&$\binom{n+i-2}i+\binom{n+i-3}{i-1}$\\[1mm]
      $a_i$  &$-\frac{i+1}q$ &$-\frac{(i+1)(w-i)(n-w-i)}{(n-2i-1)(n-2i)}$
      &$\frac{n-2+i}{n-2+2i}$\\[1mm]
   $b_i$ &$\frac {i+(q-1)(n-i)}q$& $\frac{(n+2)w(n-w)-ni(n-i+1)}{(n-2i)(n-2i+2)}$
  &$0$\\[1mm] 
   $c_i$ &$-\frac{(n-i+1)(q-1)}q$&$-\frac{(w-i+1)(n-w-i+1)(n-i+2)}{(n-2i+2)(n-2i+3)}$
  &$\frac i{n-2+2i}$\end{tabular}\vspace*{.1in}\caption{Parameters
of the metric spaces}\label{table1}
\end{table}

The following bound on $s$-codes is well known. It was proved 
by Delsarte \cite{del73,del73a} for codes in $Q$-polynomial 
association schemes which
includes $H_q^n$ and $J^{n,w}$, and
by Delsarte et al. \cite{del77b} for codes in $S^{n-1}.$ 
\begin{theorem} Let $\cC$ be an $s$-code in a compact distance-transitive 
space $X$. Then
   \begin{equation}\label{eq:bound1}
      |\cC|\le h_0+h_1+\dots+h_s.
    \end{equation}
\end{theorem}

\vspace*{.1in}For $X=S^{n-1}$ and $s=2$ this theorem gives the
bound $|\cC|\le \half n(n+3)$.
This estimate was recently improved in \cite{mus09a} where it was shown
that if the inner products
between distinct code words take values $t_1,t_2,$ and $t_1+t_2\geq 0,$
then $|\cC|\le \half n(n+1).$ The proof relies on the method of linearly
independent polynomials. Subsequently, H. Nozaki \cite{noz09a}
proved a general bound on spherical $s$-codes.
His proof builds upon Delsarte's ideas and is included here
for completeness. 

We will need a result in matrix analysis known as Ostrowski's Theorem
(\cite{hor90}, pp.224-225).
\begin{theorem} Let $F,S$ be $N\times N$ real matrices, and let $F$ be
symmetric. Let the eigenvalues of $F$ and $SS^T$ be arranged in increasing 
order, i.e., $\lambda_i(F)\le\lambda_j(F), \lambda_i(SS^T)\le\lambda_j(SS^T),
i<j.$ For each $k=1,\dots, N$ there exists a real number $\theta_k,
0\le\theta_k\le \lambda_N(SS^T)$ such that
   $$
     \lambda_k(SFS^T)=\theta_k\lambda_k(F).
   $$
\end{theorem}

\begin{theorem}\label{thm:N} \rm{(Nozaki \cite{noz09a})}
Let $\cC=\{\bfx_1,\dots,\bfx_M\}\subset X$ be an $s$-code with distances $d_1,\dots,d_s.$ 
Consider the polynomial $f(x)=\prod_{i=1}^s(d_i-x)$ and suppose that its
expansion in the basis $\{p_i\}$ has the form $f(x)=\sum_i f_ip_i(x).$ Then
    $$
      |\cC|\le \sum_{i: f_i>0} h_i.
    $$
\end{theorem}
\begin{proof} Let $|\cC|=M$ and 
consider the $M\times h_l$ matrix $H_l$ given by 
$(H_l)_{i,j}=\phi_{l,j}(\bfx_i),$ where $i=1\dots, M; j=1\dots, h_l.$
Let $\sH=(H_0,H_1,\dots,H_s)$ and consider the $M\times M$ matrix $A=\sH F\sH^t$
where 
  $$
  F=f_0I_1\oplus f_1I_{h_1}\oplus\dots\oplus f_sI_{h_s}
  $$
is a direct sum. By (\ref{eq:addition}) the
general entry of $A$ equals
$A_{\bfx,\bfy}=f(d(\bfx,\bfy)),$ which implies that
$A=f(0)I_M.$ 
\remove{Thus $A$ is positive definite. Also $M=\rank(A)\le \rank(\sH),$ and so
the rows $\bfz_i, i=1,\dots,M$ of $\sH$ are linearly independent. 
Thus we have found $M$ linearly independent vectors $\bfz_i$
that lie in the positive subspace of the quadratic form $F.$ Therefore,
their number cannot exceed
the number of positive eigenvalues of $F$ (counted with multiplicities).}

Here our arguments deviate from \cite{noz09a}.
Let $S=\big[\!\!\text{\small\begin{tabular}{c}$\sH$\\
  0\end{tabular}}\!\!\big]$ be an $N\times N$ matrix, $N=\sum_{i=0}^s h_i,$
and let $A'=SFS^T.$ The eigenvalues of $A'$ are $0$ and $f(0)$ with
multiplicities $N-M$ and $M$, respectively. 
By Ostrowski's theorem, to every positive eigenvalue of $A'$ there
corresponds a positive eigenvalue of $F,$ i.e.,
           $$
              M=|\{k: \lambda_k(A')>0\}|\le |\{k: \lambda_k(F)>0\}|,
           $$
which was to be proved.
\end{proof}

\vspace*{3mm}To apply this theorem let us compute some coefficients of the polynomial $f(x).$
\begin{lemma}\label{lemma:coeffs}
Let $f(x)=\prod_{i=1}^s (d_i-x)=\sum_{i=0}^s f_ip_i(x).$ Then
   $$
     f_s=(-1)^s r_s c_1c_2\dots c_s \quad{(s\ge 1)},
   $$
   $$
      f_{s-1}=(-1)^s r_{s-1}c_1\dots c_{s-1}\sum_{j=1}^s(b_{j-1}-d_j)
     \quad(s\ge 2).
   $$
\end{lemma}
\begin{proof}
  We have
   $$
     f(x)=(-1)^s(x^s-(d_1+\dots+d_s)x^{s-1})+\dots.
   $$
In the following we use the relation $\sL(x^m p_k)=0,$ valid for 
all $0\le m<k,$ and relations (\ref{eq:3term}) and (\ref{eq:fc}).
We compute
   \begin{align*}
     (-1)^s f_s&=r_s\sL(x^s p_s)=r_s\sL(x^{s-1} (a_sp_{s+1}+b_sp_s+c_sp_{s-1}))
       \\&=r_sc_s\sL(x^{s-1}p_{s-1})=\dots=r_sc_1c_2\dots c_s.
   \end{align*}
Next we claim that
   \begin{align*}
     \sL(x^{s}p_{s-1})=c_1c_2\dots c_{s-1}(b_0+b_1+\dots+b_{s-1}), \quad s\ge 2.
   \end{align*}
Indeed, $\sL(x^2p_1)=\sL(x (b_1p_1+c_1))=b_1c_1+b_0c_1.$ Now let us
assume that 
   $$
   \sL(x^{s-1}p_{s-2})=c_1c_2\dots c_{s-2}(b_0+b_1+\dots+b_{s-2}).
   $$
Then
   \begin{align*}
     \sL(x^{s}p_{s-1})&=\sL(x^{s-1}(b_{s-1}p_{s-1}+c_{s-1}p_{s-2}))\\
       &=b_{s-1}c_1c_2\dots c_{s-1}+
               c_{s-1}(c_1c_2\dots c_{s-2}(b_0+b_1+\dots+b_{s-2}))
   \end{align*}
as was to be proved.
Next,
   $$
    f_{s-1}=r_{s-1} \sL((-1)^s(x^s-(d_1+\dots+d_s)x^{s-1})p_{s-1})
   $$
   $$
    =(-1)^sr_{s-1} (c_1c_2\dots c_{s-1}(b_0+b_1+\dots+b_{s-1})
    -(d_1+\dots+d_s)c_1\dots c_{s-1})
   $$
   $$
    =(-1)^sr_{s-1}c_1\dots c_{s-1}((b_0+b_1+\dots+b_{s-1})-
      (d_1+\dots+d_s)).
   $$
\end{proof}

The next theorem provides an improvement of the general bound 
(\ref{eq:bound1}). It will be used in subsequent sections to establish the
main results of this paper.
\begin{theorem}\label{thm:new} Let $\cC$ be a code in a compact 
distance-transitive space $X$ with distances
$d_1,\dots, d_s.$ Let the numbers $b_i,c_i, i\ge 0$ be defined by 
(\ref{eq:3term}) and let 
   $$
      D=b_0+\dots+b_{s-1}-d_1-\dots-d_s.
   $$

{\rm (a)}. Suppose that $c_i<0, i=1,2,\dots$  and $D\ge 0.$ Then
    $$
       |\cC|\le h_0+h_1+\dots+h_{s-2}+h_s.
    $$

{\rm (b)}. Suppose that $c_i>0, i=1,2,\dots.$ Then
     $$       
        |\cC|\le \begin{cases} h_0+h_1+\dots +h_{s-2} 
                  & s\equiv1\hspace*{-3mm}\pmod 2 \text{ and } D\ge 0\\
                 h_0+h_1+\dots+h_{s-1} & s\equiv1\hspace*{-3mm}\pmod 2 
             \text{ and } D<0\\
             h_0+h_1+\dots+h_{s-2}+h_s &s\equiv0\hspace*{-3mm}\pmod 2                               \text{ and } D\le0.
         \end{cases}
      $$
\end{theorem}
\begin{proof} The proof uses Theorem \ref{thm:N} and is completed by the analysis of the signs of
$f_s$ and $f_{s-1}$ for the cases specified in the theorem.\end{proof}

{\em Remark.} It is possible to evaluate other coefficients of the
polynomial $f(x)$ in Lemma \ref{lemma:coeffs} which will lead to further
refinements of bound (\ref{eq:bound1}) from Theorem \ref{thm:N}. 
However the conditions
on the distances will involve higher-degree symmetric functions of them,
which limits somewhat their usefulness.

\vspace*{.1in}
{\em Example 1.} Consider the binary extended Golay code $\cG_{24}$ of
length $n=24$ and cardinality 4096. The distances between distinct codevectors
of $\cG_{24}$ are 8, 12, 16, and 24 \cite[p.~67]{mac91}. Since $\cG_{24}$ is a linear code,
it contains the all-zero vector and therefore also the vector
${\bf 1}=1^{24}$ of all ones. Therefore, if $\bfx$ is a codevector then
so is the vector ${\bf 1}+\bfx$. Deleting one vector from each of such pairs,
we obtain a code $\cG_{24}^{o}$ of cardinality 
$2048=\binom {24} 1+\binom {24} 3$ with distances $d_1=8,d_2=12, d_3=16.$ 
From Table 1, $b_i=12$ for all $i$, so $D=0$, and
Theorem \ref{thm:new}(a) implies that for any code $C\in H_2^{24}$ with distances
8,12,16, we have $|C|\le h_0+h_1+h_3.$ However,
    $$
(8-x)(12-x)(16-x)=\nicefrac34 K_1(x)+\nicefrac34 K_3(x),
    $$
meaning that $f_0=0,$ 
so the bound can be tightened to $|C|\le h_1+h_3.$ In other words, 
the Golay ``half-code'' $\cG_{24}^{o}$ is
an optimal $3$-distance code of length 24. This example will be generalized
in Sect.~\ref{sect:Hamming} below.

In the following sections we will use another general bound on codes 
known as Delsarte's ``linear programming'' bound \cite{del73}.
For $s$-codes this bound gives
\begin{theorem}{\rm (Delsarte)} \label{thm:Delsarte}
Let $\cC\subset X$ be an $s$-code with distances $d_1,\dots,d_s.$ 
Then 
     \begin{align*}      
|\cC|\le \max\{1+\alpha_1+\dots+\alpha_s: &\sum_{i=1}^s 
        \alpha_i p_k(d_i)\ge -p_k(\tau_0), k\ge 0;\\
           &\alpha_i\ge 0, i=1,\dots, s\}.
           \end{align*}
Here $\tau_0=0$ for the Hamming and Johnson spaces and $\tau_0=1$ for the
sphere $S^{n-1}$.
\end{theorem}

\section{Constant weight codes and intersecting families}
Call a family $\cF=\{F_1,F_2,\dots\}$ of subsets of an $n$-element set
$w$-{\em uniform} if $|F_i|=w, i=1,2,\dots,$ and call it $s$-{\em intersecting} if
$\forall_{F_i,F_j}|F_i\cap F_j|\in \{w,\ell_1,\dots,\ell_s\}$ for some
$\ell_1,\dots, \ell_s, 0\le \ell_i<w.$ 
For two subsets $F_1,F_2$ with $|F_1\cap F_2|=\ell$ the distance between their
indicator vectors $\bfx_1,\bfx_2$ equals $d_J(\bfx_1,\bfx_2)=w-\ell.$ Thus,
the indicator vectors of $\cF$ form an $s$-code $\cC$ in $J^{n,w}.$ 

\begin{theorem} \label{thm:RCW} {\rm (Ray-Chaudhuri--Wilson \cite{ray75})}
Let $\cF$ be a $w$-uniform $s$-intersecting family. Then
   \begin{equation}\label{eq:RCW}
     |\cF|\le \binom ns.
   \end{equation}
\end{theorem}
\begin{proof} 
Follows from (\ref{eq:bound1}) and Table \ref{table1}.
\end{proof}

Deza, Erd{\"o}s, and Frankl \cite{dez78} showed that for 
$n\ge w\binom {3w}w$ this estimate can be 
improved to 
 \begin{equation}\label{eq:def}
    |\cF|\le \prod_{i=1}^s \frac{n-\ell_i}{w-\ell_i}.
 \end{equation}
The particular case $\{\ell_1,\dots,\ell_s\}=
\{w-s,w-s+1,\dots,w-1\}$ corresponds to the celebrated 
Erd{\"o}s-Ko-Rado theorem \cite{erd61}. According to it, if
$n\ge (w-s+1)(s+1)$ then
  $$
  |\cF|\le \binom{n-w+s}{s}.
  $$
Note also that generally (\ref{eq:RCW})
is best possible because
the bound is met by $\cF=\binom{[n]}w.$ 
Several generalizations of Theorem \ref{thm:RCW} were obtained in 
\cite{alo91,bab95,sne94}. 
\remove{Namely, assume that 
the maximum distance in the constant weight code $\cC$ is at most $D$.
Then according to this theorem, $|\cC|\le \binom{n-w+D}{D}$ for all
$n\ge (D+1)(w-D+1).$
For instance, for $s=2,$ $d_1=1,d_2=2,$ the Erd{\"o}s-Ko-Rado bound gives 
$|\cC|\le \half(n-w+2)(n-w+1).$ This estimate coincides with part (a)
of the above corollary for $w=3$ and improves upon it for greater $w$.
On the other hand, the results in this section hold true with no conditions 
on the code's diameter.}
We obtain the following general improvement of this theorem.
\begin{theorem}\label{thm:RCWnew}
 Let $\cF$ be a $w$-uniform $s$-intersecting family.
Suppose that 
     \begin{equation}\label{eq:cond}
        \ell_1+\dots+\ell_s\ge\frac{s(w^2-(s-1)(2w-n/2))}{n-2(s-1)}.
   \end{equation}
Then
  \begin{equation}\label{eq:RCWnew}
   |\cF|\le \binom ns -\binom n{s-1}\frac{n-2s+3}{n-s+2}.
  \end{equation}
\end{theorem}
\begin{proof} The proof will follow from 
Theorem \ref{thm:new}(a). For it to hold, we need that
   \begin{equation}\label{eq:ps}
      \sum_{i=1}^s d_i=ws-\sum_{i=1}^s \ell_i\le \sum_{i=0}^{s-1} b_i.
   \end{equation}
Now take the value of $b_i$ from Table \ref{table1} and use induction
to show that 
  $$
  \sum_{i=0}^{s-1}\frac{(n+2)w(n-w)-ni(n-i+1)}{(n-2i)(n-2i+2)}=
    \frac{ws(n-w)-\binom s2 n}{n-2(s-1)}.  
  $$    
The proof is concluded by substituting this expression 
for $\sum_{i=0}^{s-1} b_i$
in (\ref{eq:ps}).
\end{proof}

Let us show that the region of $\ell_i$'s defined in (\ref{eq:cond}) is 
not void. Write this inequality as 
  $$
    \sum \ell_i> ws-\frac{ws(n-w)-\binom s2 n}{n-2(s-1)}.
  $$
As $s\le w\le n/2,$ the numerator of the fraction is nonnegative and
$n-2(s-1)\le n.$ Thus (\ref{eq:cond}) will hold if $\sum \ell_l > w^2s/n-\binom
s2.$ This last inequality holds in turn if $w$ is close to $n/2$ 
and the $\ell_i$s are large.
For instance  if $s=2$ then the Ray-Chaudhuri-Wilson
bound can be tightened for all $\ell_1+\ell_2>(2w(w-2)+n)/(n-2).$ See also
the example in the end of this section. 
The bound (\ref{eq:RCWnew}) is not as good as (\ref{eq:def}) whenever
the latter applies; on the other hand, (\ref{eq:RCWnew}) involves no
restrictions on $n$.

Let us consider in more detail the case of 2- and 3-intersecting families,
switching to the language of constant weight codes.

\begin{corollary}\label{cor:new} Let $\cC\subset J^{n,w}$ be a code.

 (a) Suppose that the distances between distinct vectors of $\cC$  
take values $d_1,d_2.$ 
 If
  \begin{equation}\label{eq:c2}
     d_1+d_2\le \frac{2w(n-w)-n}{n-2}
  \end{equation}
then
   $|\cC|\le\half (n-1)(n-2)$.

(b)  Suppose that the distances between distinct vectors in $\cC$ 
take values $d_1,d_2,d_3.$ If
  $$
     d_1+d_2+d_3\le \frac{3w(n-w)-3n}{n-4}
  $$
then    $
   |\cC|\le \frac n6(n^2-6n+11). 
    $
\end{corollary}
We note that a 2-distance constant weight code can be constructed
by taking the $\binom{n-w+2}2$ vectors
with $w-2$ ones in the first coordinates and the remaining 2 ones anywhere
outside them. This code attains the Erd{\"o}s-Ko-Rado bound and in the
case $w=3$ is extremal for Part (a) of the above corollary for all $n\ge 6.$

To establish the next result we will need the following result of Larman, Rogers, and Seidel \cite{lar77}, 
restated here in the form convenient to us: {\em Suppose that 
$\cC\subset H_2^n$ is a binary code with distances $d_1<d_2$, and $|\cC|> 2n+3.$
Then $d_1/d_2=(k-1)/k$ where
$k$ is an integer satisfying $2\le k\le \half+\sqrt{n/2}.$} 
Below we call this relation for the numbers $d_1,d_2$ the LRS condition.

\begin{proposition} \label{prop:cwc}
(a) For $6\le n\le 44$ and $3\le w\le n/2$ 
with the exception of the cases $(n,w)=(23,7),(44,17)$
the size of a 2-distance code
$\cC\subset J^{n,w}$ satisfies $|\cC|\le \half(n-1)(n-2).$

\vspace*{1mm}
(b)
If $n$ and $w$ satisfy any of the following conditions:\\[.05in]
\centerline{\begin{tabular}{cl}
  $6\le n\le 8$  &and \;$w=3$;\\
  $9\le n\le 11$ &and \;$3\le w\le 4$;\\
  $12\le n\le 14$ \;or\; $25\le n\le 34$&and\; $3\le w\le 5$;\\
  $15\le n\le 24$\; or\; $35\le n\le 46$&and\; $3\le w\le 6,$
\end{tabular}}\\

\noindent then the maximum 2-distance code 
$\cC\subset J^{n,w}$ satisfies $|\cC|=\half(n-w+1)(n-w+2).$
\end{proposition}
\begin{proof} Part (a). If the distances in $\cC$ satisfy (\ref{eq:c2}), then 
the upper bound in part (a) follows from the previous corollary.
Otherwise we examine every pair of distances $d_1,d_2$. If a given pair
does not satisfy the LRS condition, then $|\cC|\le 2n+3.$ If this condition
is satisfied, we compute the Delsarte bound of Theorem
\ref{thm:Delsarte}. Together these arguments 
produced an upper bound $\binom {n-1}2$ on the code size for all the
parameters in the statement.

Part (b). For all $n,w\le n/2$ there exists a constant weight 
2-distance code of size $\binom {n-w+2}2$. 
The matching upper estimates are established by computing the Delsarte
bound.
\end{proof}

\vspace*{-2mm}
As an example of the arguments involved in the proof, 
let $\cC$ be a two-distance code in $J^{n,w}$ with $n=13,w=5.$ There are 10 
possibilities for the distances $d_1,d_2.$ The LRS condition is fulfilled
if $d_1/d_2=(k-1)/k,2\le k\le 3.$ Thus, the
pairs $(1,3),(1,4),(1,5),$ $(2,5),(3,4),(3,5),(4,5)$ do not satisfy it,
so for all these cases $|\cC|\le 29.$ 
Next we compute the Delsarte bound $D(d_1,d_2)$ for the 3 remaining 
cases, obtaining $D(1,2)=45,D(2,3)=33,D(2,4)=27.$
This exhausts all the possible cases, so we conclude that
$|\cC|\le 45.$ As mentioned above, the extremal configuration has
45 vectors at distances 1 or 2. This code meets both
the Delsarte bound and the Erd{\"o}s-Ko-Rado bound.
This establishes both parts of the last proposition in the case considered.

Likewise, if $n=18,w=8,$ Corollary \ref{cor:new}(a) applies whenever
$d_1+d_2\le 8.$ For any such two-distance code we obtain 
$|\cC|\le 136.$
The remaining possibilities for the distances are covered by
the LRS condition or checked by computing the Delsarte bound.
This establishes the corresponding case of Part (a) of the proposition.

Generally, the Delsarte bound is better than the other bounds
for $n$ up to about 45 and is rather loose (and difficult to compute) 
for greater $n$. 

Note that the case $n=23, w=7$ is a true exception in part (a) of Prop.~\ref{prop:cwc}.
Indeed, the 253 vectors of weight 7 in the binary Golay code of length 23
have pairwise Johnson distances 4 and 6 \cite[p.~69]{mac91}, which is greater
than $\binom{23}2=231.$

\section{$s$-codes in the Hamming space}\label{sect:Hamming}
Let $\cC\subset H_q^n$ be a code in which the distances between
distinct codewords are $d_1,d_2,\dots,d_s.$
Theorem \ref{thm:new} implies the following bound.
\begin{theorem}  Suppose that 
  $$
    d_1+\dots+d_s\le \frac sq[(q-1)n-\half(q-2)(s-1)]
     \quad(\half sn \;\;\text{for $q=2$ }).
   $$ 
Then
  \begin{equation}\label{eq:H}
    |\cC|\le 1+n(q-1)+\binom n2(q-1)^2+\dots+\binom n{s-2}(q-1)^{s-2}+
    \binom ns(q-1)^s.
  \end{equation}
\end{theorem}
This enables us to draw some conclusions for sets of binary vectors with
few distances.
\begin{theorem} (a) Let $\cC$ be a binary code in which 
the distances between distinct codewords are $d_1,d_2.$ If $d_1+d_2\le n$ then 
$|\cC|\le \half(n^2-n+2).$ 

(b) Let $\cC$ be a binary code in which the distances between 
distinct codewords are $d_1,d_2,d_3.$ If $d_1+d_2+d_3\le3n/2$ 
then $$|\cC|\le 1+n+\binom n3.$$ If in addition none
or two of the three distances $d_1,d_2,d_3$ are $> n/2$ then
$$|\cC|\le n+\binom n3.$$
\end{theorem}
\begin{proof} Part (a) follows from the previous theorem. 

Part (b). Consider the annihilator polynomial 
$f(x)=(d_1-x)(d_2-x)(d_3-x)$ and let $f_0,\dots,f_3$ be its coefficients 
in the Krawtchouk basis.
We know that under the assumption of the theorem, $f_2\le0.$ This
proves the first claim in part (c). Further, by (\ref{eq:fc}), 
the constant coefficient equals
  $$f_0=-\Big(\frac n2-d_1\Big)\Big(\frac n2-d_2\Big)\Big(\frac n2-d_3\Big)
+\frac n4(d_1+d_2+d_3)-\frac{3n^2}8
  $$
If both assumptions in part (b) of the theorem hold then $
f_0\le0.$ This proves the bound $|\cC|\le n+\binom n3.$
\end{proof}

\begin{proposition}
(a) If $6\le n\le 74$ with the exception of the values $n=47,53,59,65,70,$ $71,$ or
if $n=78,$ then the size of a maximal code
with 2 distances equals $\half(n^2-n+2)$. 

(b) If $8\le n\le 22$ or $n=24$ then the size of a maximal code
with 3 distances equals $n+\binom n3$.

(c) If $10\le n\le 33$ then the size of a maximal code with 4 distances
equals $1+\binom n2+\binom n4.$ 
\end{proposition}
\begin{proof}
Part (a). Observe that the size of the code $\cC$ formed of all vectors of 
weight 2 and 
the all-zero vector equals $1+\binom n2$ for all $n\ge 3.$ It remains to show that
even if $d_1+d_2\ge n+1,$ no two-distance code of length $n$ for each value of $n$ in the statement can have
larger size. To establish this, for each $n$ 
we compute the Delsarte bound of Theorem \ref{thm:Delsarte} for all the 
possible distance values $d_1,d_2, d_1+d_2 \ge n+1$ 
that satisfy the LRS condition $d_1/d_2=(k-1)/k.$ 
These computations show that in each case the Delsarte bound 
is less than or equal to $\half(n^2-n+2)$. This establishes our claim.

Part (b). We proceed in a way analogous to part (a). Note that 
the code $\cC$ formed of all vectors of weights 1 and 3 
has size $|\cC|=n+\binom n3$ for all $n\ge 3$.  
We need to show that even if
$d_1+d_2+d_3\ge 3n/2+1,$ no three-distance code of length 
$8\le n\le 22$ or $24$ can have larger size.
To do this, we rely on Part(b) of the previous theorem. Namely, for
each $n$ in the range and for all $d_1,d_2,d_3$ such that
$d_1+d_2+d_3\ge 3n/2+1$ or that $f_0>0,$ we compute the 
Delsarte bound of Theorem \ref{thm:Delsarte} and verify that it is less
than or equal to the claimed code size.

Part (c). For $n\ge 6,$ a 3-code of size $1+\binom n2+\binom n4$ is formed of 
all vectors
of weights 0,2,4. Therefore,
if $f_1\le 0$ and $f_3\le 0$ in the expansion
  $$
  \prod_{i=1}^4(d_i-x)=\sum_{i=0}^4 f_iK_i(x),
  $$
then the claim holds true. Otherwise, for every $10\le n\le 33$ and 
for every set of numbers
$d_1,d_2,d_3,d_4$ that fails these conditions, we compute the Delsarte bound
and verify that it is less than or equal to $1+\binom n2+\binom n4.$ 
\end{proof}

Example 1 above shows that an extremal 3-distance code $\cG_{24}^{o}$ 
of length $n=24$ can be obtained from the binary Golay code $\cG_{24}.$ 
A related example accounts for the omission of $n=23$ from part (b).
Indeed, the even subcode of the Golay code $\cG_{23}$ has distances
$8,16,24,$ but its size 2048 is greater than $23+\binom {23} 3=1794,$
so this case is a true exception.

\section{Spherical codes}
Let $\cC\subset S^{n-1}$ be a code such that the inner product of any
two distinct code vectors takes one of the $s$ values
$t_1,\dots,t_s.$ Let
    $$
     M_s:=\binom{n+s-1}{s}+\binom{n+s-2}{s-1}.
    $$
\begin{theorem} {\rm (Delsarte-Goethals-Seidel, 1977)}
  $
  |\cC|\le M_s.
  $
\end{theorem}
\begin{proof} Follows from (\ref{eq:bound1}) and Table \ref{table1}
by the identity 
$\sum_{k=0}^p\binom{m+k}{k}=\binom{m+p+1}{p}$.
\end{proof}

This result was improved in \cite{mus09a} as follows:
   If $s=2$ and $t_1+t_2\ge 0$ then $|\cC|\le \half n(n+1).$
We now have the following general improvement.
\begin{theorem} Suppose that $s$ is even and $t_1+t_2+\dots+t_s\ge 0,$ then
  $$
    |\cC|\le M_{s-2}+\frac{n+2s-2}{s}\binom{n+s-3}{s-1}.
  $$
\end{theorem}
\begin{proof} Consider the polynomial $g(x)=\prod_{i=1}^s(x-t_i).$
By Lemma \ref{lemma:coeffs} its leading coefficients in the basis
of Gegenbauer polynomials are
   $$
     g_s=r_s c_1c_2\dots c_s>0, \quad g_{s-1}=r_s(-t_1-t_2\dots-t_s)\prod_i{c_i}.
   $$
Thus, $g_{s-1}\le 0$ if $t_1+\dots +t_s\ge 0$ (since $c_i>0$ for all $i$). 
Then the last case of Theorem \ref{thm:new}(b) applies, and the result follows
from Table \ref{table1}.
\end{proof}

Any binary code can be mapped to $S^{n-1}$ by a distance-preserving mapping, 
so the bound for spherical
codes implies bounds on binary codes (both constant weight and unrestricted).
However the bounds thus obtained are generally inferior to the results
derived in the corresponding discrete spaces. 
This is because the bounds become progressively stronger
as we move from a space to its subspaces, so there is no gain in
using the last theorem for binary codes.

The methods discussed in this paper are applicable to other distance-transitive
spaces of interest to geometry and combinatorics. We point to one such class
of spaces, namely, $q$-analogs of the Hamming and Johnson spaces, for
which intersection theorems were studied in \cite{fra86a,fu99}.

{\em Acknowledgment:} We are grateful to a reviewer for detailed and
insightful comments on the first version of this paper.

\def\cprime{$'$} \def\cprime{$'$} \def\cprime{$'$}
\providecommand{\bysame}{\leavevmode\hbox to3em{\hrulefill}\thinspace}
\providecommand{\MR}{\relax\ifhmode\unskip\space\fi MR }
\providecommand{\MRhref}[2]{%
  \href{http://www.ams.org/mathscinet-getitem?mr=#1}{#2}
}
\providecommand{\href}[2]{#2}

\end{document}